\documentclass[12pt]{amsart}
\usepackage{cite,amssymb,graphicx,caption,lipsum}
\usepackage{amsfonts}
\usepackage{enumerate}
\usepackage{epsfig}
\usepackage[dvipsnames]{xcolor}
\usepackage{psfrag}
\usepackage{hyperref}
\usepackage{babel}
\usepackage{lineno}
\usepackage{chemarr}
\usepackage{rotating}
\usepackage{graphicx}

\allowdisplaybreaks

\usepackage{rotating}

\usepackage[colorinlistoftodos]{todonotes}
\usepackage{pgfplots}
\usepackage{stackengine}
\usepgfplotslibrary{fillbetween}
\pgfplotsset{compat=newest}
\usetikzlibrary{decorations.markings, patterns}

\usepackage[short,nodayofweek,level,12hr]{datetime}

\usepackage{tikz}
\usetikzlibrary{math}
\usepackage{pgfplots}
\usepackage{stackengine}
\usepgfplotslibrary{fillbetween}
\pgfplotsset{compat=newest}
\usetikzlibrary{decorations.markings, patterns}

\newtheorem{prop}{Proposition}
\newtheorem{remark}{Remark}

\definecolor{Bluish}{rgb}{0.,0.,0.5}
\definecolor{Reddish}{rgb}{0.5,0.,0.}

\def\eps{\varepsilon}

\usepackage[colorinlistoftodos]{todonotes}
\usepackage[normalem]{ulem}

\definecolor{ColorTomasL}{rgb}{0.9,0.2,0.1}

\definecolor{ColorJS}{rgb}{0.9,0.,0.4}

\definecolor{ColorSanti}{rgb}{0.1,0.5,0.1}

\colorlet{Mag}{magenta!60}

\usepackage{graphicx}
\graphicspath{ {./figures/} }

\usepackage{enumerate}
\usepackage{color}
\usepackage{comment}
\usepackage{hyperref}

\newcommand{\R}{\ensuremath{\mathbb{R}}}

\renewcommand{\thefootnote}{\roman{footnote}}

\usepackage{lipsum}% http://ctan.org/pkg/lipsum
\makeatletter
\g@addto@macro{\endabstract}{\@setabstract}
\newcommand{\authorfootnotes}{\renewcommand\thefootnote{\@fnsymbol\c@footnote}}%

\makeatother

\title[WT helper virus-Defective Interfering genomes-RNA satellite dynamics]{No two without three: Modelling dynamics of the trio RNA virus-defective interfering genomes-RNA satellite}

\author[J. T. L\'azaro]{J.Tom\'as~L\'azaro\textsuperscript{1,2,3,8}}
\email{jose.tomas.lazaro@upc.edu}
\author[A. Alb\'o]{Ariadna Alb\'o\textsuperscript{1}}
\author[T. Alarc\'on]{Tom\'as Alarc\'on\textsuperscript{4,2,5,8}}
\author[S.F. Elena]{Santiago F. Elena\textsuperscript{6,7}}
\author[J. Sardany\'es]{Josep Sardany\'es\textsuperscript{2,8}}
\email{jsardanyes@ crm.cat}

\address{\textsuperscript{1}Departament de Matemàtiques, Universitat Politècnica de Catalunya, Avda. Diagonal 647, 08028 Barcelona, Spain
\\
\textsuperscript{2}Centre de Recerca Matem\`atica (CRM). Campus de Bellaterra. Edifici C, Cerdanyola del Vallès 08193 Barcelona, Spain
\\
\textsuperscript{3}Institute of Mathematics of the UPC-BarcelonaTech (IMTech), C. Pau Gargallo 14, 08028 Barcelona, Spain
\\
\textsuperscript{4}ICREA, Pg. Llu\'is Companys 23, 08010 Barcelona, Spain
\\
\textsuperscript{5}Departament de Matem\`atiques, Universitat Aut\`onoma de Barcelona. Campus de Bellaterra. Edifici C, Cerdanyola del Vall\`es 08193 Barcelona, Spain
\\
\textsuperscript{6}Instituto de Biolog\'ia Integrativa de Sistemas (I$^2$SysBio), CSIC-Universitat de Val\`encia, Parc Cient\'ific UV, Avda. Catedr\'atico Agust\'in Escardino 9, Paterna, 46980 Val\`encia, Spain
\\
\textsuperscript{7}The Santa Fe Institute, 1399 Hyde Park Road, Santa Fe, NM 87501, USA
\\
\textsuperscript{8}Dynamical Systems and Computational Virology, I$^2$SysBio-CRM Associated Unit}

\begin{document}

\maketitle

\begin{abstract}
Almost all viruses, regardless of their genomic material, produce defective viral genomes (DVG) as an unavoidable byproduct of their error-prone replication. Defective interfering (DI) elements are a subgroup of DVGs that have been shown to interfere with the replication of the wild-type (WT) virus. Along with DIs, other genetic elements known as satellite RNAs (satRNAs), that show no genetic relatedness with the WT virus, can co-infect cells with WT helper viruses and take advantage of viral proteins for their own benefit. These satRNAs have effects that range from reduced symptom severity to enhanced virulence. The interference dynamics of DIs over WT viruses has been thoroughly modelled at within-cell, within-host, and population levels. However, nothing is known about the dynamics resulting from the nonlinear interactions between WT viruses and DIs in the presence of satellites, a process that is frequently seen in plant RNA viruses and in biomedically relevant pathosystems like hepatitis B virus and its $\delta$ satellite. Here, we look into a phenomenological mathematical model that describes how a WT virus replicates and produces DIs in presence of a satRNA at the intra-host level. The WT virus is subject to mechanisms of  complementation, competition, and various levels of interference from DIs and the satRNA. Examining the dynamics analytically and numerically reveals three possible stable states: (i) full extinction, (ii) satellite extinction and virus-DIs coexistence and (iii) full coexistence. Assuming DIs replicate faster than the satRNA owed to their smaller size drives to scenario (ii), which implies that DIs could wipe out the satRNA. In addition, a small region of the parameter space exists wherein the system is bistable (either scenarios (ii) or (iii) are concurrently stable). We have identified transcritical bifurcations in the transitions between scenarios (i) to (iii) and saddle-node bifurcations behind the change from bistability to monostability. Despite the model simplicity, our findings may have applications in biomedicine and agronomy. They will cast light on the dynamics of this three-species system and aid in the identification of scenarios in which the clearance of the satRNAs may be possible thus \emph{e.g.}, allowing for less severe disease symptoms.

\vskip1cm


\end{abstract}

\section{Introduction}

Viruses are found infecting organisms from all realms of the Tree of Life. Viruses are obligate intracellular parasites that lack of translation machinery to complete their infection cycles. Hence, they need to infect and take profit of the cell's machinery to replicate their genomes and produce the structural proteins that will be used for packaging their genomes.  Perhaps the most remarkable characteristic of viruses, in particular those having RNA genomes, is their high mutation rate, consequence of a lack of proof-reading mechanisms in their replicases \cite{Sanjuan2010a}. At the one hand, this high mutation rate, along with their very short generation time and large population size, bestow viral populations with great evolvability \cite{Andino2015}.  At the other hand, RNA viruses' extremely compacted genome organization makes mutations potentially harmful. In fact, most randomly introduced mutations either impose a significant fitness cost or are fatal \cite{Sanjuan2010b}. These highly deleterious or lethal mutations can vary from point mutations to genomic deletions of variable length; these mutants are collectively referred as defective viral genomes (DVGs) \cite{Vignuzzi2019}.

A fraction of deletion DVGs has been long shown to interfere with genome replication and accumulation, being known as defective interfering (DI) RNAs. DI RNAs were first reported by Preben von Magnus~\cite{VonMagnus1954}, who studied their accumulation in influenza A virus populations passaged in embryonated chicken eggs. Based on these serial passage experiments the existence of incomplete virus variants which increase rapidly in frequency and cause drops in overall virus titers was proposed. The existence of virus variants with large genomic deletions has been confirmed thereafter in many virus families \cite{Vignuzzi2019}, both with RNA and DNA genomes. DI RNAs are thought to replicate much faster than full-length wild-type (WT) viruses, due to their smaller genome sizes. Moreover, DI RNAs can evolve other strategies to better compete with WT viruses. DI RNAs cannot autonomously replicate because they lack most, if not all, of WT coding sequences. They must, therefore, co-infect a cell with a WT virus in order to replicate, becoming obligate parasites of WT viruses. As the frequency of the DIs increases, the overall virus production is reduced because essential WT-encoded gene products are no longer available (\emph{i.e.}, interference) \cite{Chao2017}. DI RNAs can have implications for virus amplification in cultured cells, protein expression using viral vectors, and vaccine development \cite{PalmaHuang1974}. Nearly all animal and many plant RNA viruses infections are associated to DI RNAs. The viral genes necessary for movement, replication, and encapsidation are typically absent from these truncated and frequently rearranged versions of WT viruses, but they still have all of the \emph{cis}-acting components needed for replication by the WT virus's RNA-dependent RNA polymerase (RdRp). 

In the past 20 years, \emph{de novo} generation of DI RNAs has received a great deal of attention. For plant virus DI RNAs, the RdRp-mediated copy choice model, which was first outlined for the generation of DI RNAs from animal viruses, still holds true \cite{White1999}. DI RNAs are probably subject to intense selective pressure for biological success after \emph{de novo} generation. While the majority of DI RNAs attenuate the WT virus's symptoms, DI RNAs of broad bean mottle virus and of turnip crinkle virus (TCV) possess the unusual attribute of exacerbating symptom severity (reviewed in \cite{Simon2004} and in \cite{Badar2021}). Interestingly, DI RNAs can also be produced by DNA viruses such as hepatitis B virus (HBV). Defective forms of HBV, named spliced HBV, have been characterized and investigated \emph{in vivo} \cite{TerrePetit1991,RosmorducPetit1995}. HBV DNA genome is transcribed into a pre-genomic RNA (pgRNA) by the viral P protein in the cell nucleus. Then pgRNAs are exported to the cytoplasm to be further processed to produce mature viruses. During the synthesis of pgRNA molecules, P also produces defective RNAs \cite{TerrePetit1991} which, after reverse transcription in the cytoplasm result in defective DNA genomes, can be packaged into mature viral particles, thus behaving as DI agents.

Another relevant member of the subviral RNA brotherhood are the so-called satellite viruses and the satellite RNAs (satRNAs), which can be either linear or circular (also known as virusoids) \cite{Simon2004,Palukaitis2016}. While satellite viruses generally encode for the components to build their own capsid protein, but depend on the helper WT virus for replication and movement, satRNAs often do not encode for any protein. Typically, virus satellites have been suggested to establish symbiotic relations with the WT helper virus, thus getting a benefit. However, other side-effect processes such as competition may arise during co-infection. Moreover, some satellite viruses can also act as parasites of the WT virus, thus taking a profit of the presence of the WT virus but not providing an advantage to it.  SatRNA and satellite virus genomes are mostly or completely unrelated to their WT helper virus genome, a major difference with DI RNAs.  The diversity of satRNAs and satellite virus structure and interaction with their helper WT virus is remarkable. For example, satC associated with TCV is a hybrid molecule composed of sequence from a second satRNA and two portions from the $3'$ end of TCV genomic RNA~\cite{Simon2004}. The satRNA associated to the ground rosette virus (GRV) further confounds earlier classifications. While not necessary for viral movement within a host, this noncoding satRNA is necessary for GRV to encapsidate in the coat protein of its luteovirus partner, groundnut rosette assistor virus, as a requirement for aphid transmission~\cite{Robinson1999}. Although more often found in plant viruses, some satellites are known to infect vertebrates~\cite{Krupovic2016}, insects~\cite{Ribiere2010}, and unicellular eukaryotic cells~\cite{Schmitt2002}. Some virus satellites have a strong clinical impact.  For example, HBV has its own satellite RNA virus, the hepatitis $\delta$ virus (HDV). Infections with HBV are more virulent, quickly evolving towards fatal cirrhosis when there is coinfection with HDV~\cite{Taylor2020}. HDV is replicated by cellular RNA polymerases I, II and III but uses for packaging HBV envelope proteins in order to accomplish viral particle assembly and release~\cite{Taylor2020}.
 
Understanding the host's reaction to viral invasion has recently made strides that have helped to clarify how DI RNAs, satRNAs and satellite viruses cause, enhance, or minimize disease symptoms. For instance, symptom attenuation was once primarily ascribed to direct competition for limited replication factors between the helper WT virus and subviral RNA \cite{Chao2017}. Recent data from a number of viral species, however, point to the possibility that the enhancement of host resistance by subviral RNA may be just as important, if not more so \cite{Palukaitis2016,Gnanasekaran2019,Badar2021}. Concepts defining the genetic connection between WT viruses and subviral RNA are also developing. A recent study suggests that some pairs of subviral RNA and helper WT viruses have more complex relationships, including mutualistic ones benefiting both participants \cite{Simon2004,Badar2021}. Furthermore, in natural infections, WT viruses could support the replication of more than one subviral element.  For example, these three-ways interactions shall be relevant to understanding the dynamics of HBV, HBV-derived DI RNAs and HDV. Even more complex systems exist, as it is the case for panicum mosaic virus (PMV), which is found coinfecting with a satellite virus (sPMV) and at least two satRNAs (S and C) and DI RNAs produced both from the WT virus as well as from sPMV \cite{Qiu2001,Pyle2018}, or the case of the bipartite tomato black ring virus that coexists with DIs derived from its both genomic RNAs as well as with a satRNA that affects its vertical transmission efficiency; all the interactions being strongly dependent on the host species \cite{Hasiow2018,Pospieszny2020,Minicka2022}.
\begin{figure}[h]
\begin{centering}
\captionsetup{width=\linewidth}
\includegraphics[width=0.89\linewidth]{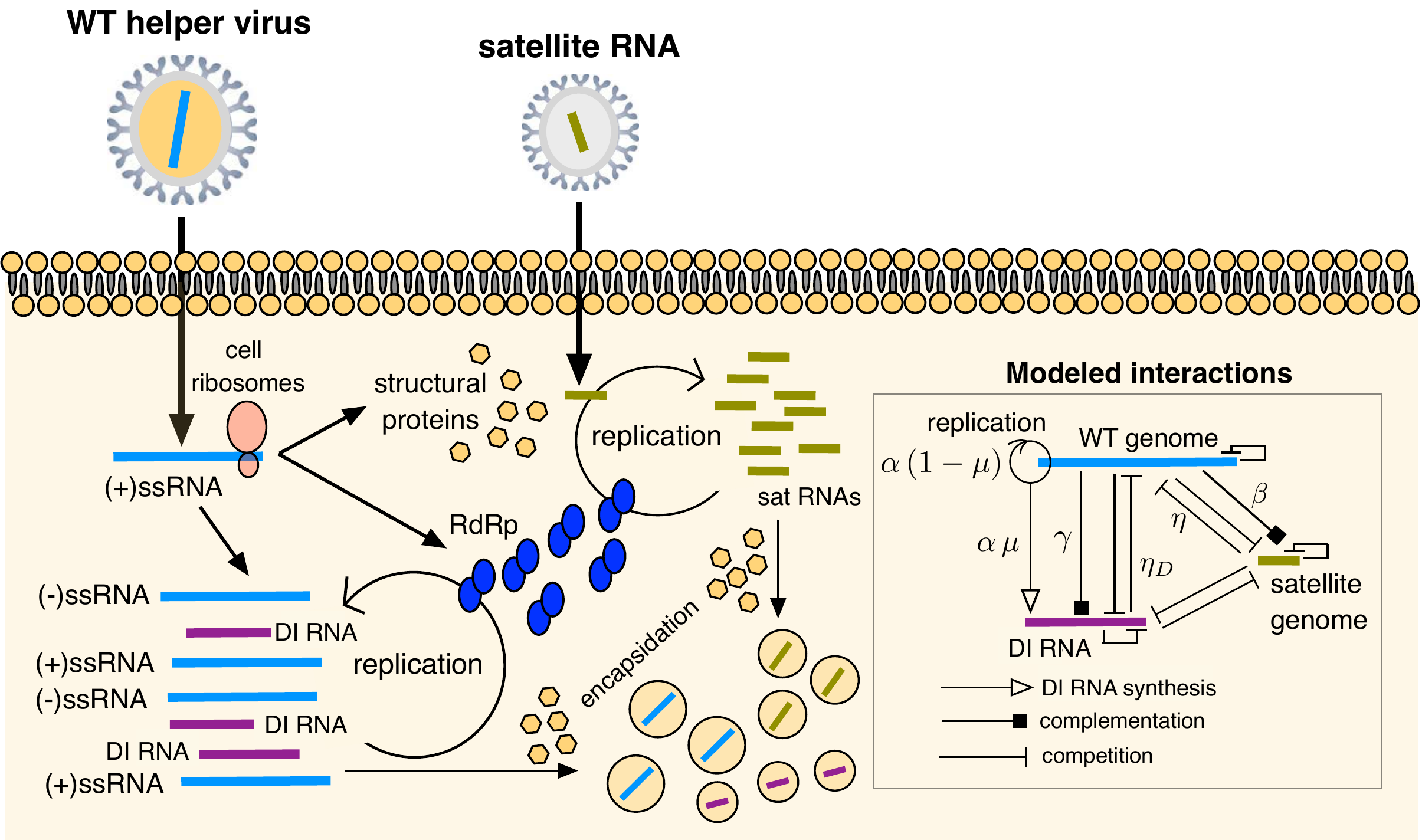}
\caption{{\small Schematic diagram of the interactions between a positive-sense single-stranded RNA virus (wild type helper virus, labeled HV) that produces DI RNAs during its replication and co-infects with a satRNA.  The inset displays the interactions considered in Eqs.~\eqref{eq1}-\eqref{eq3}, including synthesis of DI RNAs at a rate $\mu$, complementation of DI RNAs and the satRNA with the products synthesized by the HV (RNA-dependent RNA-polymerase (RdRp) and coat proteins, not considered explicitly), and competition between all RNA types. The terms $\eta_D$ and $\eta$ are included to study different interference strengths exerted by the DI RNAs and the satRNA on the HV.}} \label{fig:model}
\end{centering}
\end{figure}

Theoretical investigations of the dynamical impact of DI RNAs in the replication of WT viruses have been carried out by several authors~\cite{Szathmary1993,BanghamKirkwood1990,Kirkwood1994,Sardanyes2010,Chao2017}. Typically, these mathematical models had taken mean-field approximations considering either discrete-~\cite{Szathmary1993,Zwart2013} or continuous-time \cite{BanghamKirkwood1990,Kirkwood1994,Chao2017} dynamical systems. However, to the extend of our knowledge, the only previous theoretical study incorporating both helper and satellite viruses used an epidemiological approach in which host individuals could be infected by different combinations of viral and subviral RNAs~\cite{LuciaSanz2022}. Here we take a population dynamics approach to explore the within-host dynamics of a system of molecular replicators composed by a WT helper virus, one satRNA and the DI RNAs generated during WT virus replication (Fig.~\ref{fig:model}). With this approach, we want to determine the dynamics arising from most basic principles of replication and interaction between replicators without entering into mechanistic details involving proteins. By doing so, satRNAs and viral satellites could be considered as homologous.  For simplicity, hereafter we will refer to the WT helper virus as HV. 

The manuscript is organised as follows. In Section~\ref{se:model} we introduce the mathematical model. Section~\ref{se:anal:results} contains analytical results concerning the domain of the dynamics, equilibrium points and their local stability. Section~\ref{se:numerical:analysis} illustrates different scenarios from numerical results, also including information on transients for those systems with DI RNAs shorter than the satRNA and for which no full coexistence is possible. Finally, we show the system can display bistability and that achieving either satRNA clearance with HV-DIs persistence or full coexistence depends on the initial populations of replicators.

\section{Mathematical model}
\label{se:model}

We develop a dynamical model based on three coupled autonomous ordinary differential equations (ODEs) to investigate the dynamics of a WT helper virus (HV) population supporting the replication of a satRNA together with the synthesis of DIs as a by product of the replication of the HV genome. Let us denote by $x=(V,S,D)$ the following state variables being the (normalised) concentration of HV ($V$), the satRNA ($S$) and, for simplicity, all possible DI RNAs grouped into a single category ($D$), respectively. The corresponding system of ODEs is given by:
\vspace{0.1cm}
\begin{eqnarray}
\dot{V} &=& \alpha \, (1-\mu) \,V \, \Omega(x) - \eps V, \label{eq1} \\
\dot{S} &=& \beta \, V \, S \, \theta(x) - \eps \, S, \label{eq2} \\
\dot{D} &=& \left( \alpha \, \mu  + \gamma  \, D \right)\, V \,  \theta(x) - \eps \, D, \label{eq3}
\end{eqnarray}
with $\Omega(x) = 1-V-\eta S - \eta_D D,$ and $\theta(x) = 1-V-S-D$.

We will refer in short to the model as $\dot{x}=F(x)$. The model considers well-mixed populations and takes into account the processes of virus replication, complementation, competition with asymmetric interference strengths, and spontaneous degradation of the different RNAs (see Fig.~\ref{fig:model} for a schematic diagram of the modeled processes). To keep the model as simple as possible, the production of viral proteins is ignored and replication/encapsidation processes for the satRNA and the DIs are made proportional to the amount of HV (simulating complementation). The replication rates of the viral genomes are proportional to parameters $\alpha$ (HV), $\beta$ (satRNA), and $\gamma$ (DIs). We will generically assume that $\beta, \gamma > \alpha$. This assumption is based on the fact that both DIs and the satRNA genomes are always shorter than the genome of the HV (see tables 1 - 5 in \cite{Badar2021}), and thus replication is expected to be faster. For example, tobacco mosaic virus has a genome size of ca. 6.4 kb and its satellite virus sTMV is about 1.1 kb~\cite{Arenal1999};  TCV genome is 4.1 kb long while its satC has only 0.4 kb~\cite{Altenbach1981}. Lucerne transient streak virus is about 4.2 Kb long while its satellite scLTSV is 0.3 Kb long~\cite{Badar2021}; HBV pgRNA is about 3.5 kb long, while HDV is 1.7 kb. Interestingly, in the case of HBV, the length of the DI RNA (deletion-containing pgRNA) is about 2.2 kb~\cite{Rosmorduc1995}, a bit longer than the satellite HDV. As we will show below, the case where DIs replicate faster that the satRNAs ($\gamma > \beta$) does not allow for the coexistence of the three populations. For those cases with satRNAS replicating faster than DIs  ($\beta>\gamma$) coexistence is possible. This latter case may correspond to viruses supporting very short satRNAs such as linear or circular ones.

The replication of the HV unavoidably results in the production of DIs at a rate $\mu$ (we assume $0< \mu <1$). Both $\Omega(x)$ and $\theta(x)$ are logistic functions introducing competition between the three viral populations due to finite host resources. Notice that the logistic function for the HV, given by $\Omega(x)$, involves the competition parameters $\eta_D, \eta >1$ to investigate higher interference strengths by the satRNAs and the DIs on the HV. Such interference may be due to competition for host resources, viral components shared by the three RNAs (\emph{e.g.}, envelop proteins) or triggering of host antiviral defenses by an excessive accumulation of viral particles or post-transcriptional gene silencing in response to the accumulation of different RNA species \cite{Palukaitis2016,Gnanasekaran2019,Badar2021}. Finally, parameter $\eps$ denotes the degradation rate of all RNA molecules, which, for simplicity, is considered to be the same for the three populations considering that the expected growth asymmetries have been introduced in replication rates.

\section{Analytical results}
\label{se:anal:results}

In this section we first study the domain where dynamics take place and compute the nullclines of the system. Then, we provide an analysis of its equilibrium points~\eqref{eq1}-\eqref{eq3} and their local stability. %We illustrate the identified dynamics of the model by tuning some parameter values, focusing on the production rate of DIs ($\mu$) and on the interference coefficients of both DIs and the satRNAs over the HV.

\subsection{Domain of confined dynamics and nullclines}
\label{se:domain:nullclines}
As it is common in many biological models, the competition for limited resources term $\theta(x)$ limits populations' growth and confines the dynamics to a finite domain. In our case, this is given by the tetrahedron
\begin{equation}
\mathcal{U} = \left\{
x=(V,S,D) \ \bigg| \ x\geq 0 \quad \textrm{and} \quad V+S+D \leq 1 
\right\},
\label{def:domain:U}
\end{equation}
which is determined by the coordinate planes and the plane $\theta(x)=0$ \emph{i.e.}, the plane $V+S+D=1$. The fact that the planes $V=0$ and $S=0$ are invariant under the dynamics of~\eqref{eq1}-~\eqref{eq3} and that the vector field $F$ of~\eqref{eq1}-\eqref{eq3} points inwards in the rest of its faces, makes the domain $\mathcal{U}$ positively invariant. That is, orbits with initial conditions on $\mathcal{U}$ remain inside of this domain for all $t\geq 0$.

Let us now compute the nullclines of Eqs.~\eqref{eq1}-~\eqref{eq3}, which determine the regions of increase or decrease of the variables. In our case, the nullcline $\dot{V}=0$ is easily computable and exhibits two connected components: the planes $V=0$ and
\begin{equation}
\Omega(x) = \frac{\eps}{\alpha (1-\mu)} 
\Leftrightarrow
V+ \eta S + \eta_D D = \sigma,
\label{V:nullcline}
\end{equation}
where $\sigma$ is defined as
\begin{equation}
\sigma:= 1 - \frac{\eps}{\alpha (1-\mu)}.
\label{def:sigma}
\end{equation}
The $V$-nullcline component $V=0$ is biologically trivial: the absence of HV leads to no satRNA and no DIs (since both need the first to be present) and, therefore, to total extinction.
\begin{figure}
\label{fig:Udomain}
\begin{centering}
\captionsetup{width=\linewidth}
\includegraphics[width=0.45\linewidth]{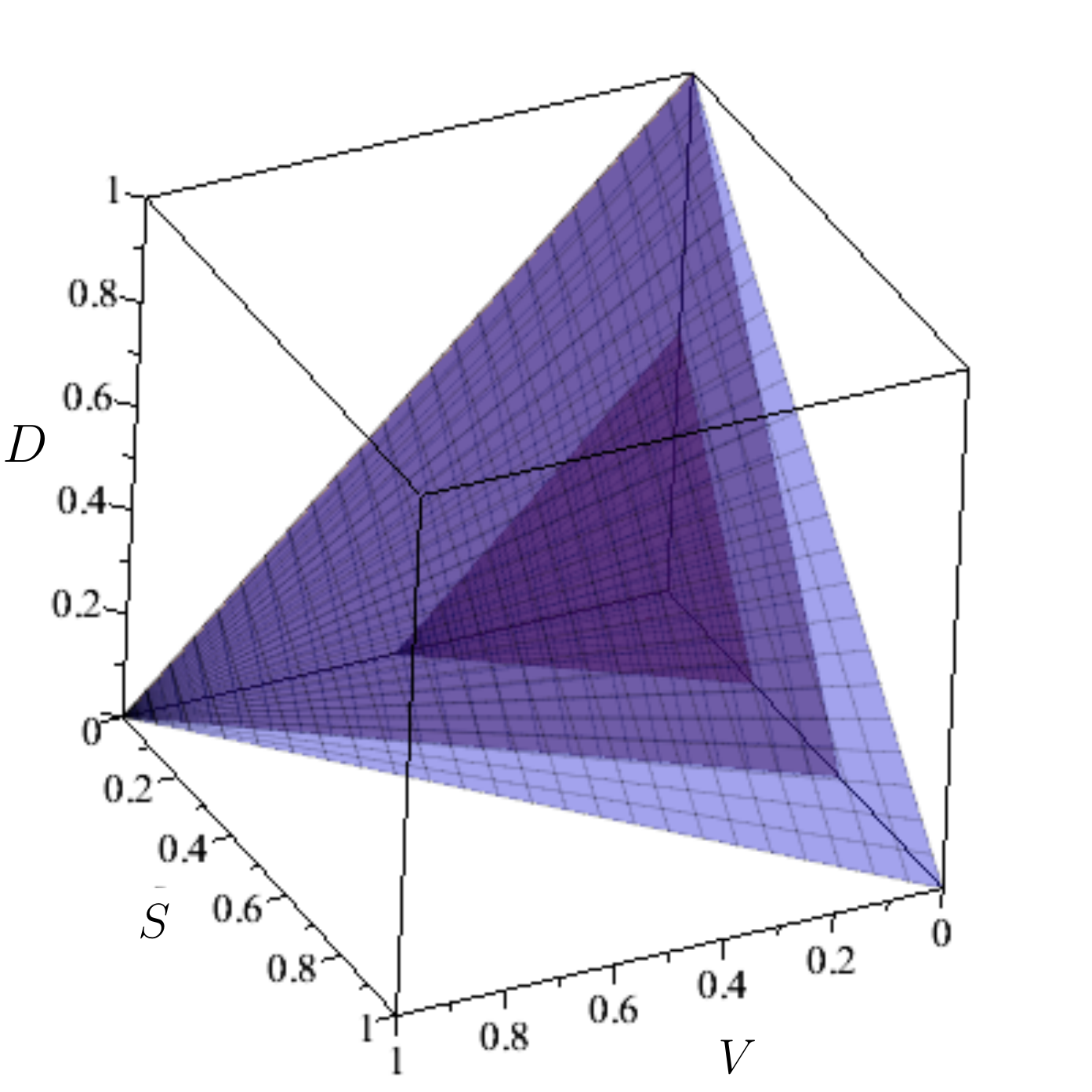}
\caption{
{\small From outside to inside (also, from lighter to darker): the planes
$V+S+D=1$ (\emph{i.e.} $\theta(x)=0$) , 
$V+\eta S+ \eta_D D=1$ (\emph{i.e.} $\Omega(x)=0$) , and
$V+\eta S+ \eta_D D= \sigma$ (\emph{i.e.} $\Omega(x)=1-\sigma$) are depicted. The parameter $\sigma$ is defined in~\eqref{def:sigma}.}}
\end{centering}
\end{figure}
The second one~\eqref{V:nullcline} determines the evolution of the WT helper virus in $\mathcal{U}$. A view of this domain and the latter planes are depicted in Fig.~\ref{def:domain:U}. The other two nullclines, $\dot{S}=0$ and $\dot{D}=0$ do not provide a simple representation. The first one is formed by the (invariant) plane $S=0$ and the (piece of) hyperbolic cylinder
\[
\mathcal{U} \cap \left\{ V \, (1-V-S-D) = \frac{\eps}{\beta} \right\}. 
\]
The second one, $\dot{D}=0$, is given by the algebraic surface
$\left( \alpha \mu + \gamma D \right) V \, \theta(x) = \eps D$. 
Notice also that $\Omega(x) < \theta(x)$ for any $x=(V,S,D)\in \R_+^3\setminus{(1,0,0)}$ and that $\Omega(1,0,0) = \theta(1,0,0)$.

\subsection{Equilibrium points and local stability}
\label{se:equilibrium-points}
The equilibrium points are the solutions $x^*$ of $F(x)=0$. As usual, their (local) stability is approached through its linearized system $\dot{x}=DF(x^*)(x-x^*)$, whose Jacobian matrix is given by
\begin{equation}
DF(x) = \left(
\begin{array}{ccc}
\alpha(1-\mu) (\Omega(x) - V)-\eps & -\alpha (1-\mu)\eta V &
-\alpha(1-\mu)\eta_D V \\
\beta S (\theta(x) - V) & \beta V(\theta(x)-S) - \eps & - \beta VS \\
(\alpha \mu+\gamma D) (\theta(x)-V) & -(\alpha \mu + \gamma D)V & \gamma V\theta(x) - (\alpha \mu + \gamma D) V - \eps
\end{array}
\right).
\label{DF}
\end{equation}
Regarding the equilibrium points of system~\eqref{eq1}--\eqref{eq3}, the following statements hold:
\begin{itemize}
\item The origin, $\mathcal{O}=(0, 0, 0)$, is always an equilibrium point for any value of the parameters. It represents the full extinction of $V$, $S$ and $D$. Its Jacobian matrix
\[
DF(0,0,0)=\left(
\begin{array}{crr}
\alpha (1-\mu) - \eps & 0 & 0 \\
0 & -\eps & 0 \\
\alpha \mu & 0 & - \eps
\end{array}
\right)
\]
has eigenvalues $\lambda_1=\alpha(1-\mu)-\eps$, $\lambda_2=\lambda_3=-\eps<0$ (semisimple). Notice that its stability depends on the sign of $\lambda_1$. Precisely, $\lambda_1<0$ (\emph{i.e.}, $\mathcal{O}$ is locally attractor) is equivalent to the condition
\begin{equation}
\mu >  \mu_c = 1-\frac{\eps}{\alpha} \Leftrightarrow \frac{\eps}{\alpha(1-\mu)}> 1 \Leftrightarrow \sigma<0.
\label{cond:origin:gas}
\end{equation}
If this condition holds, that is, the DIs generation rate $\mu$ exceeds the critical value $\mu_c$, then all the points in $\mathcal{U}$ satisfy $\dot{V}<0$. This, in its turn, leads to $\dot{S}<0$ and $\dot{D}<0$ and, afterwards, to total extinction. Consequently, the origin is the unique equilibrium point of system~\eqref{eq1}-\eqref{eq3}, and it is a global asymptotically attractor.

\bigskip

Henceforth, let us assume that condition
\begin{equation}
\mu < \mu_c\Leftrightarrow 0 < \frac{\eps}{\alpha(1-\mu)}<1 \Leftrightarrow \sigma>0,
\label{origin:saddle:cond}
\end{equation}
is satisfied. Hence, the origin is a saddle equilibrium point, with a $2$-dimensional stable manifold and a $1$-dimensional unstable curve. The latter is tangent to the vector $(1-\mu,0,\mu)$ at the origin. In this case we also have:

\item No equilibrium points of the form $(0,S,D)$. As previously mentioned, $V=0$ leads necessarily to total extinction.

\item No equilibrium points on the line $\{ S=0, D=0 \}$ except the origin. Indeed, if this was the case they would be solutions of
\begin{eqnarray*}
\dot{V} &=& \alpha (1-\mu) V (1-V) - \eps V = 0, \\
\dot{D} &=& \alpha \mu V (1-V) = 0.
\end{eqnarray*}
From the latter equation we get either $V=0$ or $V=1$. The first case corresponds to extinction. The second one, $V=1$, leads to $\eps=0$, which does not hold since $\eps>0$ by assumption.

\item In a similar way it can be proved that there are no equilibrium points in the plane $\{ D= 0 \}$ other than the origin. Indeed, in this plane the system becomes
\begin{eqnarray*}
\alpha(1-\mu)V (1-V-\eta S)-\eps V   &=& 0, \\
\beta VS (1-V-S) - \eps S &=& 0, \\
\alpha \mu V (1-V-S) &=& 0.
\end{eqnarray*}
Since $V\neq 0$, from the last equation it turns out that $V+S=1$. Substituting into the second equation we get $S=0$ (since $\eps\neq 0$) and, therefore, $V=1$. Clearly, $(1,0,0)$ does not satisfy the first equation.

\end{itemize}
The next two propositions summarise the type of non-trivial equilibrium points that the system can have.

\begin{prop}[\textsf{No-satRNA equilibria, $P$-point}]
\label{prop:P:points}
Let us assume that condition~\eqref{origin:saddle:cond} holds. 
Then, there exists a unique biologically meaningful equilibrium point $P=(V_1,0,D_1)$ of the  system~\eqref{eq1}--\eqref{eq3}. This $P$-point satisfies that
\[
V_1=\sigma - \eta_D D_1,
\]
where $D_1$ is the unique real root in the interval $(0,\frac{\sigma}{\eta_D})$ of the following polynomial of degree $3$:
\[
q(D)=-\gamma \eta_{D} \left( \eta_{D}-1 \right) {D}^{3}+A_2 D^2+ A_1  D+\alpha\mu\sigma \left( 1-\sigma \right),
\]
with
\begin{eqnarray*}
A_2 &=& \left( 
-\alpha \mu \eta_D+\gamma\,\sigma \right)  \left( \eta_D-1 \right) -\gamma\eta_D \left( 1-\sigma \right)  \\
A_1 &=& \alpha\mu\sigma \left( \eta_D-1 \right) + \left( - \alpha\mu\eta_D+\gamma\sigma \right)  \left( 1-\sigma \right)  -\eps.
\end{eqnarray*}
In particular, this point $P=(V_1,0,D_1)$ does not depend on the parameter $\eta$.
\end{prop}

\begin{proof}
The plane $\{ S=0\}$ (absence of satRNA) is invariant under the dynamics of system~\eqref{eq1}--\eqref{eq3}. These dynamics are governed by equations
\begin{eqnarray}
\dot{V} &=& \alpha (1-\mu) V (1-V-\eta_D D) - \eps V, \label{eq:s0:1}\\
\dot{D} &=& (\alpha \mu + \gamma D) V (1-V-D) - \eps D. \label{eq:s0:2}
\end{eqnarray}
Thus, $P$-points correspond to the solutions making these equations vanish. Since $V_1>0$, the first one becomes $V_1+\eta_D D_1= \sigma$  which, in particular, implies that $0< D_1 < \frac{\sigma}{\eta_D}$.

Substituting $V_1+\eta_D D_1= \sigma$ into equation
$(\alpha \mu + \gamma D) V (1-V-D) - \eps D=0$ 
it turns out that $D_1$ must to be a root of the polynomial
\[
q(D)=(\alpha \mu + \gamma D)(\sigma - \eta_D D) \big(1-\sigma - (\eta_D-1)D \big) - \eps D
\]
(in the interval $(0,\frac{\sigma}{\eta_D})$). On one hand, since $q(0)=\alpha \mu \sigma (1-\sigma)>0$ (recall that $0<\sigma<1$) and $q(\sigma/\eta_D)=-\eps \sigma/\eta_D<0$, we get from Bolzano's theorem the existence of, at least, one zero of $q(D)$ in this interval. On the other, expanding and collecting in powers of $D$, we reach the following equivalent expression for $q(D)$:
\[
-\gamma \eta_{D} \left( \eta_{D}-1 \right) {D}^{3}+A_2 D^2+ A_1  D+\alpha\mu\sigma \left( 1-\sigma \right) 
\]
where
\begin{eqnarray*}
A_2 &=& \left( 
-\alpha \mu \eta_D+\gamma\,\sigma \right)  \left( \eta_D-1 \right) -\gamma\eta_D \left( 1-\sigma \right),  \\
A_1 &=& \alpha\mu\sigma \left( \eta_D-1 \right) + \left( - \alpha\mu\eta_D+\gamma\sigma \right)  \left( 1-\sigma \right)  -\eps.
\end{eqnarray*}
From the fact that $\eta_D>1$, it follows that $\lim_{D\to +\infty} q(D)=-\infty$ and that $\lim_{D\to -\infty} q(D)=+\infty$. So Bolzano's theorem ensures that the three roots of $q(D)$ are real: one is negative, a second one stays in the interval $(0,\frac{\sigma}{\eta_D})$ and the third one is greater than $\frac{\sigma}{\eta_D}$. Consequently, since $V_1=\sigma-\eta_D D_1 \in (0,1)$ if $D_1 \in (0, \frac{\sigma}{\eta_D})$, there is exactly one biologically meaningful $P$-point $(V_1,0,D_1)$.
\end{proof}

\begin{prop}[\textsf{Coexistence equilibrium points, $Q$-points}] \label{prop:Q:points}
Let us assume condition~\eqref{origin:saddle:cond} is satisfied. Then, $Q=(V_2,S_2,D_2)$ is a coexistence equilibrium point of system~\eqref{eq1}--\eqref{eq3} ($Q$-point in short) if and only if $Q\in \mathcal{U}$ and the following conditions hold:

\begin{itemize}
\item[(i)] its $D$-component is given by
\[
D_2=\frac{\alpha \mu}{\beta - \gamma}
\]
which, necessarily, implies that $\beta>\gamma$. In order to make $Q$, in principle, biologically meaningful, it must satisfy necessarily that
$D_2< \frac{\sigma}{\eta_D}.$

\item[(ii)] The component $0<V_2<\sigma$ is a root of the degree-$2$ polynomial
$V_2^2 + M V_2 + m=0$
where
\begin{equation}
M:=\frac{\sigma - (\eta_D - \eta) D_2 - \eta}{\eta-1}, \qquad
m:= \frac{\eta\eps}{\beta(\eta-1)}>0.
\label{def:M:m}
\end{equation}
More precisely, $V_2$ is given by
\begin{equation}
V_2^{\pm} = \frac{-M \pm \sqrt{M^2 - 4m}}{2}, 
\label{Qpoints:eq:V2}
\end{equation}
provided that $M^2-4m\geq 0$. 

\medskip

\item[(iii)] The component $0<S_2<1$ is given by the expression
\begin{equation}
S_2 = 1-V_2-D_2 - \frac{\eps}{\beta V_2},
\label{def:s2}
\end{equation}
where $0<V_2<1$ is a solution of $V_2^2 + M V_2 + m=0$.
\end{itemize}
\end{prop}

\begin{remark}
(a) The restriction $D_2<\frac{\sigma}{\eta_D}$ follows from the same argument used for the $P$-points: since the equilibrium point must fall onto the $V$-nullcline $V + \eta S + \eta_D D =\sigma$, $D$ cannot exceed this value.
(b) From statement (i) it turns out that if the DIs replication rate $\gamma$ is larger than the satRNA's, $\beta$,  coexistence $Q$-equilibria no longer exist (indeed, $D_2<0$).
(c) The maximal number of biologically meaningful $Q$-points for fixed values of the parameters is $2$. As it will be showed in the numerics, there are examples with none, one and two $Q$-points.
\end{remark}

\bigskip

\begin{proof}
\begin{itemize}
\item[(i)] We seek for points of type $Q=(V_2,S_2,D_2)\in \mathcal{U}$, with $V_2, S_2, D_2>0$, steady state of our system~\eqref{eq1}--\eqref{eq3}. In particular, this implies that $Q$ must belong to the intersection of the nullclines which are not coordinate planes. That is, $Q$ must satisfy the following three conditions:
\begin{equation}
\Omega(Q)=\frac{\eps}{\alpha(1-\mu)}, \qquad 
V_2\, \theta(Q)=\frac{\eps}{\beta}, \qquad \textrm{and} \qquad 
(\alpha \mu + \gamma D_2) V_2 \,\theta(Q)=\eps D_2.
\label{eq:conditions:coexistence:eq}
\end{equation}
Substituting the second equation into the third one it leads to
\[
(\alpha \mu + \gamma D_2) \frac{\eps}{\beta}=\eps D_2 \Rightarrow
(\alpha \mu + \gamma D_2) = \beta D_2
\]
and therefore
\begin{equation}
D_2 = \frac{\alpha \mu}{\beta - \gamma}.
\label{eq:D2}
\end{equation}
where $\beta>\gamma$ to have $D_2>0$. Since $V,S,\eta,\eta_D$ are all positive, and 
$V_2+\eta S_2 + \eta_D D_2 =\sigma$  it turns out that
$\eta_D D_2 < \sigma \Rightarrow D_2 < \frac{\sigma}{\eta_D}$.

\item[(ii)] Consider now the two first conditions in~\eqref{eq:conditions:coexistence:eq} and the value $D=D_2$ in~\eqref{eq:D2}:
\begin{eqnarray}
\Omega(Q)=1-\sigma &\Rightarrow& V_2+\eta S_2 = \sigma - \eta_D D_2  \label{eq:cond:coex1}\\
V_2 \, \theta(Q)=\frac{\eps}{\beta} &\Rightarrow&
V_2 + S_2 = 1- D_2 - \frac{\eps}{\beta v_2}. \label{eq:cond:coex2}
\end{eqnarray}
Subtracting~\eqref{eq:cond:coex1} multiplied by $\eta$ to~\eqref{eq:cond:coex2}, and performing some trivial algebraic manipulations, it turns out that $V_2$ must be a root of the following degree $2$ polynomial 
$V_2^2 + M V_2 + m =0$,
where
\begin{equation*}
M=\frac{\sigma - (\eta_D - \eta) D_2 - \eta}{\eta-1}, \qquad
m= \frac{\eta\eps}{\beta(\eta-1)}>0.
\end{equation*}
That is, $V_2$ is given by
\begin{equation*}
V_2^{\pm} = \frac{-M \pm \sqrt{M^2 - 4m}}{2}, 
\end{equation*}
provided that $M^2-4m\geq 0$. 

\item[(iii)] Once determined $V_2$ we seek an expression for $S_2$. Indeed,
\[
V_2 \theta(Q)=\frac{\eps}{\beta} \Longleftrightarrow 1-V_2-S_2-D_2 = \frac{\eps}{\beta V_2}, 
\]
and so,
\[
S_2 = 1 - V_2 - D_2 - \frac{\eps}{\beta V_2}. 
\]
The points $Q=(V_2,S_2,D_2)$ obtained in this way will be biologically meaningful provided that $Q\in \mathcal{U}$.
\end{itemize}
\end{proof}
The complex dependence (in the sense of the number of parameters involved) of the expressions of the $P$ and $Q$-points makes cumbersome to analytically determine their regions of existence and their local stability. In the next section we perform a numerical study of these equilibrium points for particular choices of the parameters. We have focused on, under view, which are the most virologically-relevant parameters (production of DI RNAS $\mu$, and interference coefficients $\eta,\eta_D$). We believe with this choice of parameters we are illustrating the most remarkable features in terms of asymptotic and transient dynamics, and bifurcation phenomena.

%%%%
\begin{figure}[ht]
\begin{centering}
\captionsetup{width=\linewidth}
\includegraphics[width=\linewidth]{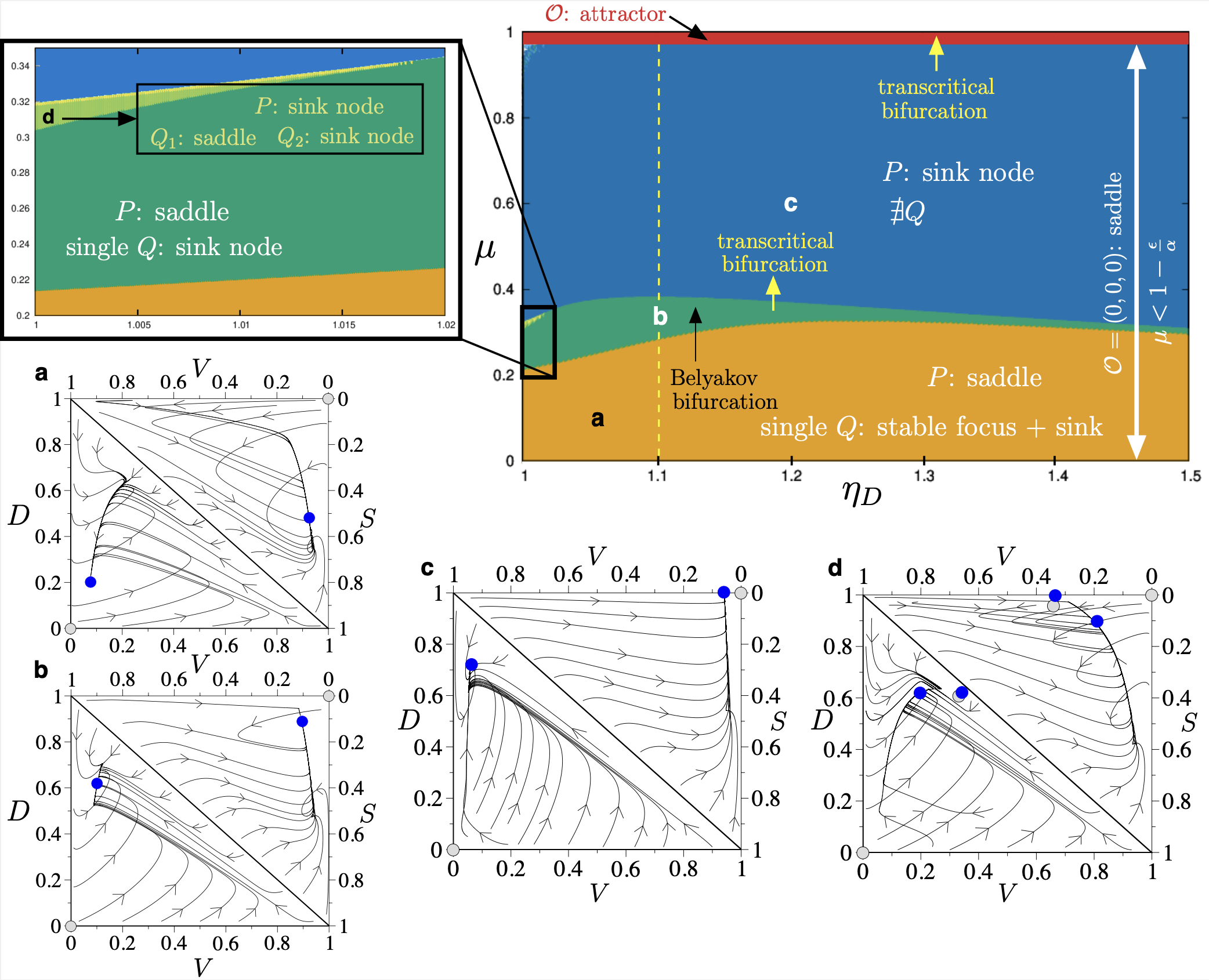}
\caption{
{\small (Main panel, top-right) Regions of existence of equilibria and stability in the parameter space $(\eta_D, \mu)$.  The inset shows a very narrow region with bistability (light green) containing equilibria $Q_1$ (saddle point) and $Q_2$ (stable focus + sink), together with $P$ (sink node). The thick white arrow indicates the region where the origin is unstable ($\mu < \mu_c = 1 - \eps/\alpha$). Figure~\ref{fig:vaps:PQ:etaD:1p1} shows information about equlibria and stability of $P$ and $Q$ points along the yellow dashed vertical line at $\eta_D = 1.1$. Eight phase portraits are shown below for the following $\eta_D$ and $\mu$ values (indicated with the same letters in the $(\eta_D, \mu)$ space): (a) $\eta_D=1.05$, $\mu=0.1$;  (b) $\eta_D=1.1$, $\mu=0.32$;  (c) $\eta_D=1.2$, $\mu=0.6$;  and (d) $\eta_D=1.001$, $\mu=0.31$. Notice in (d) that the two $Q$ points appear very close each other. Here the arrows indicate the directions of the orbits and blue and gray dots denote stable and unstable equilibria, respectively.} %In all of the plots we have used $\alpha=1$, $\beta=2$, $\eta=1.3$, $\eps=3\cdot 10^{-2}$, and $\gamma=1.5$.}
} \label{fig:Ppoints:existence:stability}
\end{centering}
\end{figure}

\section{Numerical results}
\label{se:numerical:analysis}
Numerical integration has been done with the 7th-8th order Runge-Kutta-Fehlberg-Simó method with automatic step size control and local relative tolerance $10^{-15}$. In most of the numerical results we will use initial conditions $(V(0), S(0), D(0))=(0.1, 0.05, 0)$. These initial conditions seem feasible in terms of real virus populations: an initial small quantity of HV, a lower order of magnitude quantity of its satRNA and no DIs at all that will be produced during HV replication. Despite this choice, we must note that in most of the identified scenarios initial conditions are not really important since the system is monostable. In the small region of bistability we have identified (see below) the basin of attraction of $P$-point is extremely small.

\subsection{Analysis of \boldmath{$P$-} and $Q$-points in terms of $\mu$ and $\eta_D$}
\label{se:numerics:PQpoints}
This section is devoted to the study of the equilibrium points of the system~\eqref{eq1}-\eqref{eq3}, assuming all the parameters fixed
except $\mu$ (DIs generation rate during imperfect replication of the HV) and $\eta_D$ (interference competition strength exerted by DI RNAs on the HV). The other parameters are set, if not otherwise specified, as follows:
\begin{equation}
\alpha=1, \quad \beta=2, \quad \eta=1.3, \quad \eps=3\cdot 10^{-2},
\quad \gamma=1.5.
\label{numerics:fixed:parameters}
\end{equation}
We will let $\mu \in [0,1]$ and $\eta_D \in (1,1.5]$.
Notice that, in this particular case,
\[
\mu_c = 1-\frac{\eps}{\alpha} =0.97.
\]

\begin{figure}[t]
\begin{centering}
\captionsetup{width=\linewidth}
\includegraphics[width=\textwidth]{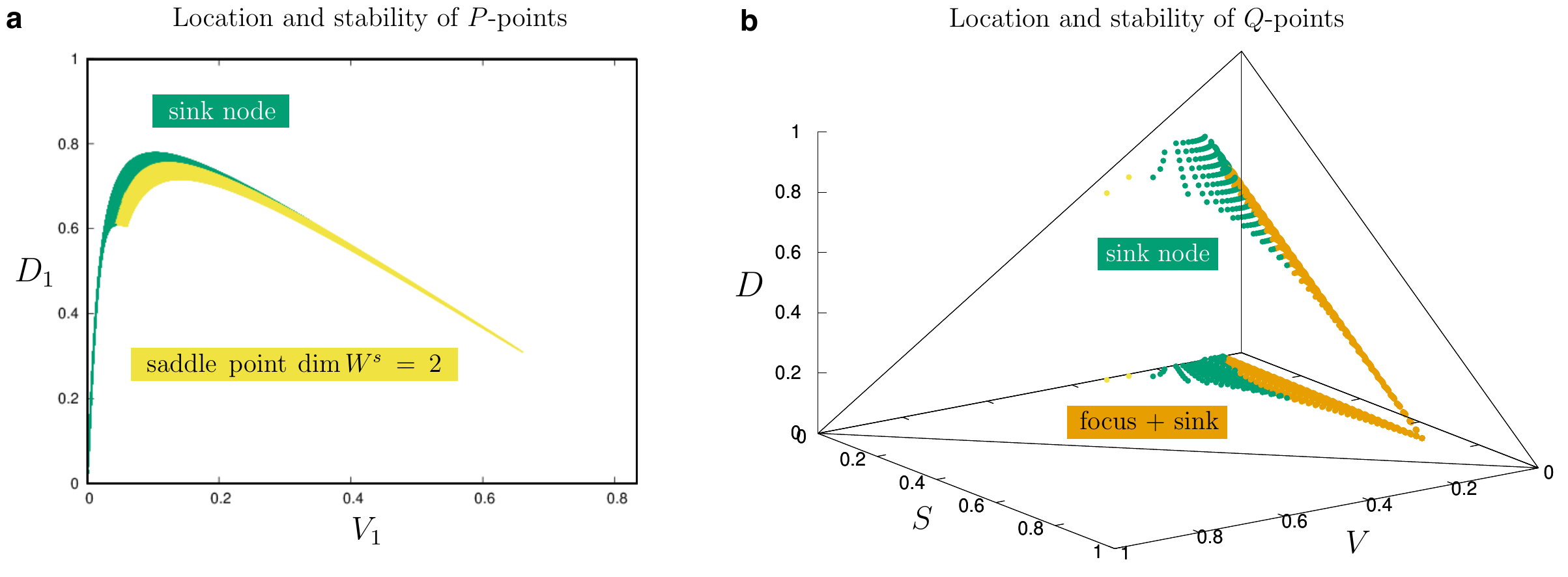}
\caption{
{\small (a) Location and stability of the $P$-point on the phase plane $(V_1, D_1)$. (b) Location and stability of the the $Q$-points in the domain $\mathcal{U}$. At the bottom, their projection on the $(V,S)$-plane. Attractors are shown in green (node) and orange (focus + sink). In yellow we display unstable points.}} \label{fig:location:PQpoints}
\end{centering}
\end{figure}

\begin{figure}
\begin{centering}
\captionsetup{width=\linewidth}
\includegraphics[width=0.9\textwidth]{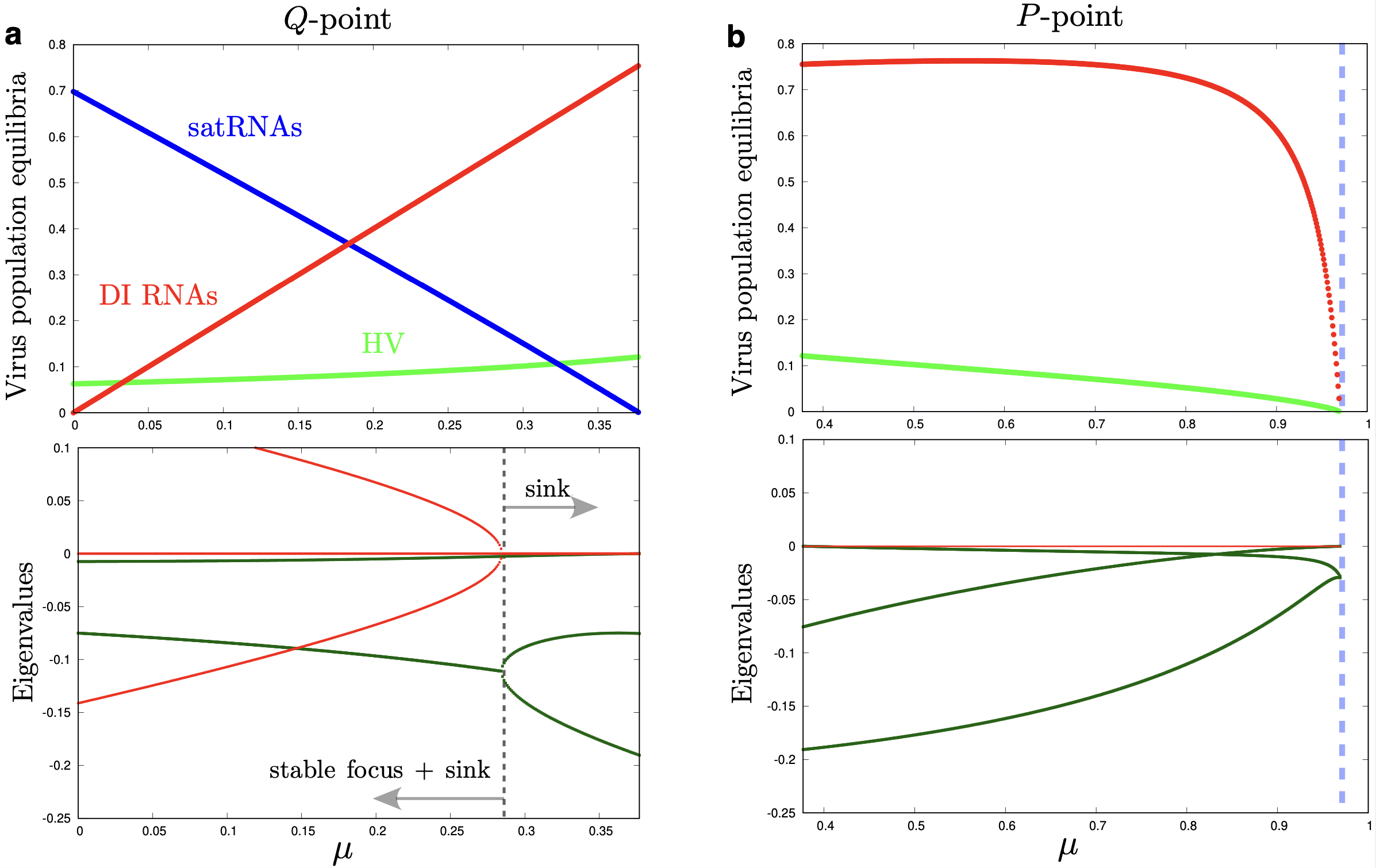}
\caption{
{\small (Top) Plot of the coordinates $(V_2,S_2,D_2)$ (a) and $(V_1,0,D_1)$ (b) of the unique attractors with $\eta_D=1.1$ for $\mu \in [0,1]$ (see vertical dashed yellow line in Fig.~\ref{fig:Ppoints:existence:stability}). We show \textcolor{green}{WT helper virus (HV)}, \textcolor{red}{DI RNAs}, and \textcolor{blue}{satRNA}. In (a) $\mu \in[0,\mu_*]$, $\mu_*\simeq 0.377$ being the value where the points $Q$ and $P$ collide in a transcritical bifurcation. (b) Plot of the $P$-point coordinates, becoming attracting after the transcritical bifurcation. This point crosses the origin at $\mu = 0.97$, in another transcritical bifurcation (thick blue dashed line). (Bottom) Eigenvalues of the Jacobian matrix at the corresponding attracting point ($Q$ left, $P$ right). Real parts of these eingenvalues are drawn in green and imaginary parts in red. For $\mu \in [0,\mu_{**})$, with $\mu_{**}\simeq 0.28$, two of them are complex (conjugate) with negative real part and a negative real eigenvalue (stable focus + sink). For $\mu \in (\mu_{**},\mu_*]$ all three eigenvalues are real negative, that is, the $Q$-point is a sink (the dashed gray line indicates the Belyakov bifurcation). Analogously for the $P$-point, all eigenvalues are real negative thus $P$ is a sink (attractor).
%Similarly to the precedent figures, the rest of parameters have been fixed to $\alpha=1$, $\beta=2$, $\eta=1.3$, $\eps=3\cdot 10^{-2}$, and $\gamma=1.5$.
}} \label{fig:vaps:PQ:etaD:1p1}
\end{centering}
\end{figure}

As already mentioned in Section~\ref{se:equilibrium-points}, if $\mu>0.97$ then the origin, the total extinction of $V$, $S$, and $D$, is a global attractor. So let us assume, henceforth, that $\mu < 0.97.$ The study we provide here is divided into several parts: existence, location in phase space and stability of equilibrium points and their bifurcations. Results on times ``to equilibrium'' (understood in the sense ``up to a given distance from it'') have  been deferred to Section~\ref{se:time:to:equilibria}, specially for those cases with outcompetition of satRNAs by the HV and DI RNAs. As mentioned, this is an interesting scenario from a biomedical point of view, especially for those systems in which the clearance of the satRNA may avoid most severe disease outcomes.

\begin{figure}
\begin{centering}
\captionsetup{width=\linewidth}
\includegraphics[width=\textwidth]{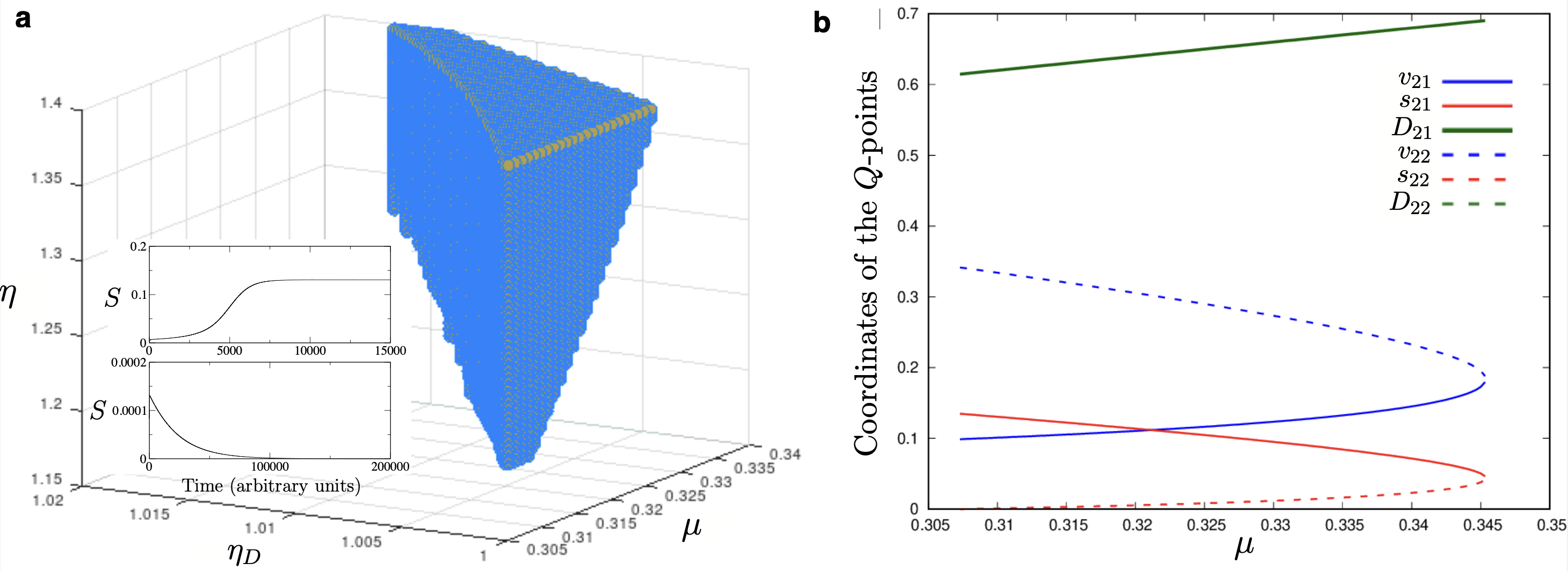}
\caption{
{\small (a) Small region in the parameter space with bistability (blue volume). The inset displays time series for the satRNA with two different initial conditions: one achieving equilibrium $Q$ (upper) and the other reaching equilibrium $P$ (bottom). (b) Bifurcation diagram displaying the coordinates of the equilibrium points and the transition from bistability to monostability at increasing $\mu$. Here $Q_{1} =(V_{21}, S_{21}, D_{21})$ is a node (continuous lines) and $Q_{2}=(V_{22}, S_{22}, D_{22})$ is a saddle (dashed line). Equilibrium point $P$ is not shown. We have used $\eta_D=1.001$, and $\eta=1.8.$}} \label{fig:bistability}
\end{centering}
\end{figure}

From Proposition~\ref{prop:P:points} we have the existence of a unique $P$-point for all $(\eta_D,\mu) \in (1,1.5] \times [0, 0.97)$. Similarly, from Proposition~\ref{prop:Q:points} we get the existence of $Q$-points in some areas inside. Precisely, in the regions depicted in \textcolor{orange}{orange} and \textcolor{ForestGreen}{green} in Figure~\ref{fig:Ppoints:existence:stability} we have a unique $Q$-point. Moreover, in a tiny region on the left hand side of the same figure, depicted in light-green, we have two coexisting $Q$-points, named $Q_1$ and $Q_2$ (see also the inset). Fixed the value of $\eta_D$, as we increase $\mu$ the (unique) $Q$-point existing in the orange and green regions approaches the plane $S=0$ and leaves the (biologically meaningful) domain $\mathcal{U}$ undergoing a collision with the corresponding $P$-point. That is, as the production rate of DI RNAs increases, the system evolves to a behaviour with a unique equilibrium point $P$ with no satRNA. As we will see below, this no-satellite equilibrium is point attractor.

Regarding the local stability of the $P$ and $Q$-points, we refer the reader to Fig.~\ref{fig:Ppoints:existence:stability} for the colors' meaning. For the sake of illustration, 
Fig.~\ref{fig:vaps:PQ:etaD:1p1} shows equilibria for the HV, satRNAs and DI RNAs, and the eigenvalues at the equilibrium points $P$ and $Q$ for the particular case $\eta_D = 1.1$ (see vertical dashed yellow line in Fig.~\ref{fig:Ppoints:existence:stability}). More generically, we have:

\begin{itemize}
\item In the \textcolor{orange}{orange} region, $P$ is a saddle point (so unstable) with a $1$-dimensional stable manifold, \emph{i.e.}, $\dim W^s(P) = 1$. There is also a unique coexistence equilibrium $Q$, which is an attractor. In particular, of type stable focus $+$ sink (so its Jacobian matrix having a couple of complex eigenvalues with negative real part and a third one real negative).
\item In the \textcolor{ForestGreen}{green} region the $P$-point is a saddle (with $\dim W^s(P) = 1$) and the $Q$-point is a sink (all three eigenvalues of its Jacobian matrix are real negative), attractor.

\item The separation between the \textcolor{orange}{orange} and the \textcolor{ForestGreen}{green} regions is given by a so-called \emph{Belyakov bifurcation curve}. This kind of bifurcation corresponds to $Q$ passing from stable focus $+$ sink to a sink. That is, the two complex eigenvalues with negative real part become real (and negative). It does not imply any change in its local stability.

\item Between the \textcolor{ForestGreen}{green} and the \textcolor{RoyalBlue}{blue} regions, $P$ (saddle) and $Q$ (sink) undergo a transcritical bifurcation: they collide and exchange their stability. Figure~\ref{fig:location:PQpoints} displays their spatial location in the plane $S=0$ (for the $P$-points) and in $\mathcal{U}$ (for the $Q$-points).

\item Inside the \textcolor{RoyalBlue}{blue} area the $Q$ point has some negative component, and so it is out of the biologically meaningful domain $\mathcal{U}$. The remaining unique equilibrium is of $P$-type, so with no-satellite, and is an attracting sink.

\item As above mentioned, in the \textcolor{Red}{red} area the origin (total extinction) is a global attractor. In the rest of the diagram, it always exists as an equilibrium but is of saddle type, with $\dim W^s(\mathcal{O})=2$.

\item Finally, in the narrow \textcolor{green}{light-green} area (see the inset in the main panel of Fig.~\ref{fig:Ppoints:existence:stability}), the system exhibits coexistence of two attracting equilibrium points: two attracting sinks points of type $P$ and $Q$, and a second unstable $Q$-point (of saddle type). In spite of its small measure, we provide some notes on this bistability scenario below.
\end{itemize}

\begin{figure}
\begin{centering}
\captionsetup{width=\linewidth}
\includegraphics[width=0.675\linewidth]{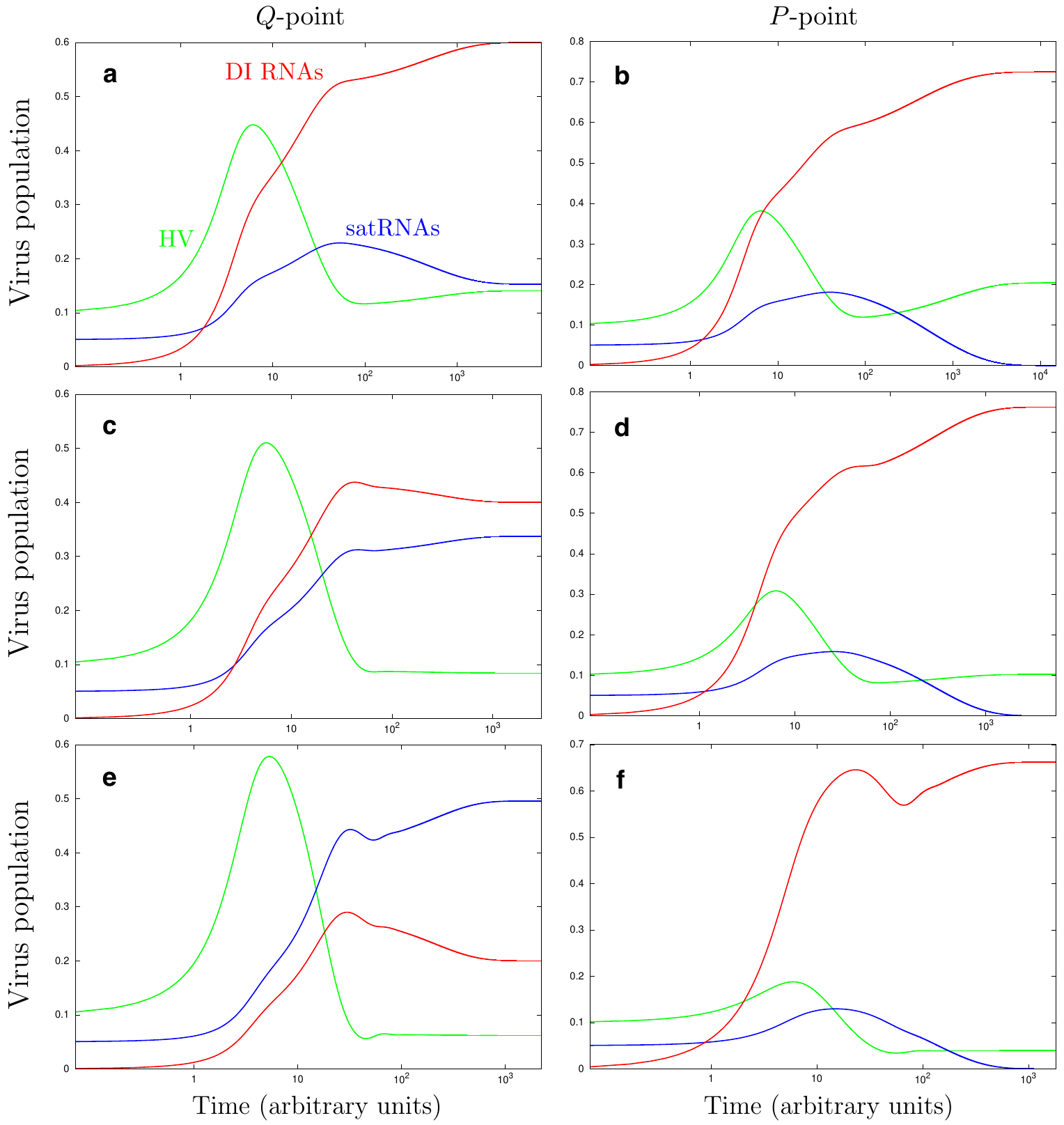}
\caption{
{\small Time series for the orbit starting with initial conditions $(V(0),S(0),D(0))=(0.1, 0.05, 0)$ and up to reaching a distance of $10^{-6}$ of the corresponding attractor. Several cases are displayed for parameter values inside the phase diagram displayed in Fig.~\ref{fig:Ppoints:existence:stability}: $\eta_D=1.03$ for $\mu=0.4$ [(a), with $Q$ sink node] and $\mu=0.3$ [(b) with $P$ sink node]; $\eta_D=1.1$ for $\mu=0.2$ [(c), $Q$ being a stable focus + sink)] and  $\mu=0.5$ [(d), with $P$ sink node)]; $\eta_D=1.3$ $\mu=0.1$ [(e), $Q$ being a stable focus + sink] for $\mu=0.7$ [(f), with $P$ sink node).% The rest of the parameters have been fixed to $\alpha=1$, $\beta=2$, $\eta=1.3$, $\eps=3\cdot 10^{-2}$, and $\gamma=1.5$.
}} \label{fig:timeseries:P:Q}
\end{centering}
\end{figure}

A very important feature of biological dynamical systems is whether they simultaneously have different stable states. This means that, for given parameter values, different initial conditions can drive to different equilibrium values (which have their own basins of attraction). Typically, biological systems can be monostable or bistable. Whether a system is bistable or monostable can have deep implications in the nature of the bifurcations and, in virology, it can involve the clearance of a given population, since often some of the two possible stable states has some component equal to zero. 

As we already mentioned, we have identified a very narrow region in the parameter space where bistability is found (Fig.~\ref{fig:Ppoints:existence:stability}). This region has been explored in more detail and is displayed in Fig.~\ref{fig:bistability}(a) by plotting those parameter values in the space $(\mu,\eta_D,\eta)$ giving place to bistability. The inset in panel (a) shows two time series using two different initial conditions. The upper one reaches the $Q$-point while the lower one the $P$-point. Despite this result indicating that systems with HV, DI RNAs, and RNA satellites could be bistable, we have to notice that the basin of attraction of the $P$-point is extremely small (results not shown), and thus the clearance of the satRNA could take place when the amount of co-infecting satRNA is very low. The transition from this bistable scenario to the region where the $P$-point is a sink node (blue area in Fig.~\ref{fig:Ppoints:existence:stability}) is governed by a saddle-node bifurcation between the $Q_1$ and $Q_2$ points (Fig.~\ref{fig:bistability}(b)). However, in the results displayed in Fig.~\ref{fig:Ppoints:existence:stability} the bifurcation curve separating coexistence from satRNA extinction is mainly due to transcritical bifurcations.

\subsection{Time to equilibria.}
\label{se:time:to:equilibria}
In order to complement the numerical analyses displayed so far, we show some characteristic time series for parameter values falling inside the phase diagram of Fig.~\ref{fig:Ppoints:existence:stability}. Notice that in all of the cases displayed in Fig.~\ref{fig:timeseries:P:Q} the population of HV starts increasing rapidly, being followed by increase in the populations of both DI RNAs and satRNAs. Once these two latest populations achieve large numbers, the HV population starts decreasing due to the their interference. The DI RNAs asymptotically achieve large population numbers and satRNAs as well at increasing their interferent effect, as shown in panels (a), (c), and (d) where $\eta_D$ has been increased. Notice that satRNAs population in panel (e) become dominant due to the large value of $\eta_D$ and the low production of DI RNAs ($\mu = 0.1$). Generically, the increase in satRNAs population seems to have a higher impact on the population of DI RNAs than on the HV. Increasing the production of DI RNAs typically involves the clearance of the satRNAs, as shown in panels (b), (d), and (f) in Fig.~\ref{fig:timeseries:P:Q}. This increase involves larger DI RNAs amounts and thus the outcompetition of the satRNAS. 

\begin{figure}
\begin{centering}
\captionsetup{width=\linewidth}
\includegraphics[width=\textwidth]{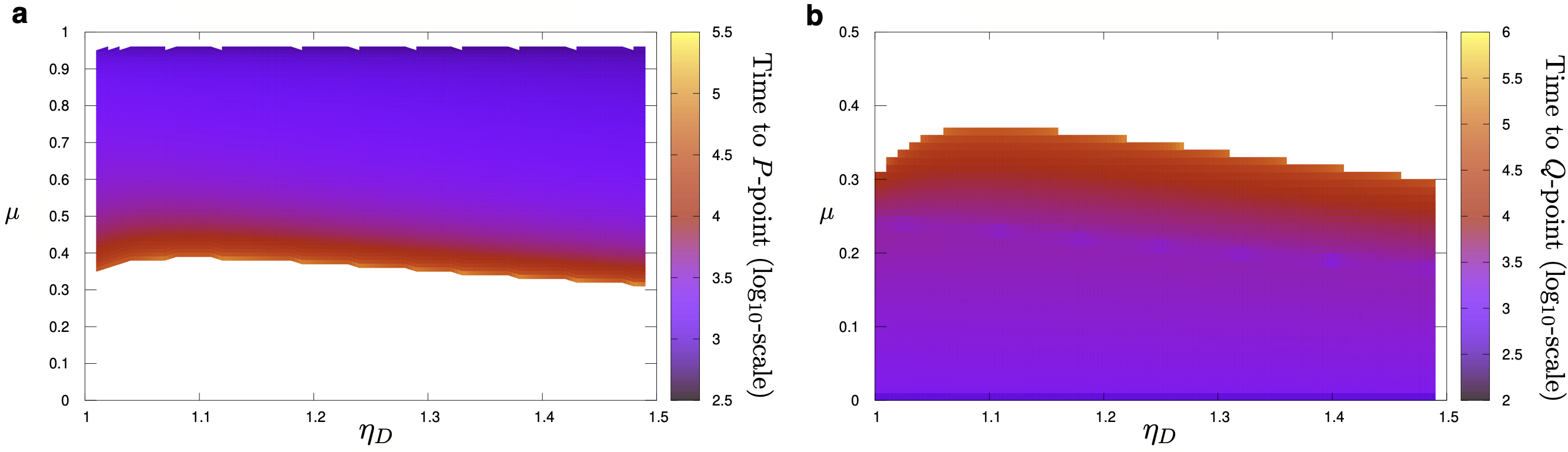}
\caption{
{\small
Transient times (in $\log_{10}$-scale) of the orbit starting with initial conditions $(V(0),S(0),D(0))=(0.1, 0.05, 0)$ to reach a distance of $10^{-6}$ of the corresponding attracting $P$-point (left) or $Q$-point (right).
The values of the parameters are the same as in the previous figures.}} \label{fig:timetoQpoints}
\end{centering}
\end{figure}

Concerning transients, as expected, they become longer close to the bifurcations separating the full coexistence from the scenario where only HV and DI RNAs persist (Fig.~\ref{fig:timetoQpoints}. This means that the time that satRNAs are able to persist in the scenario where they are outcompeted by the HV and DI RNAs depends on parameter values. For instance, Fig.~\ref{fig:timetoQpoints}(a) shows that for large values of DIs production the times are very fast (about $10^3$ time units). However, when $\mu$ is close to the transcritical bifurcation curve such transients can be about two orders of magnitude longer. The same occurs for the scenario with full coexistence found at further  decreasing $\mu$.

\begin{figure}
\captionsetup{width=\linewidth}
\includegraphics[width=\textwidth]{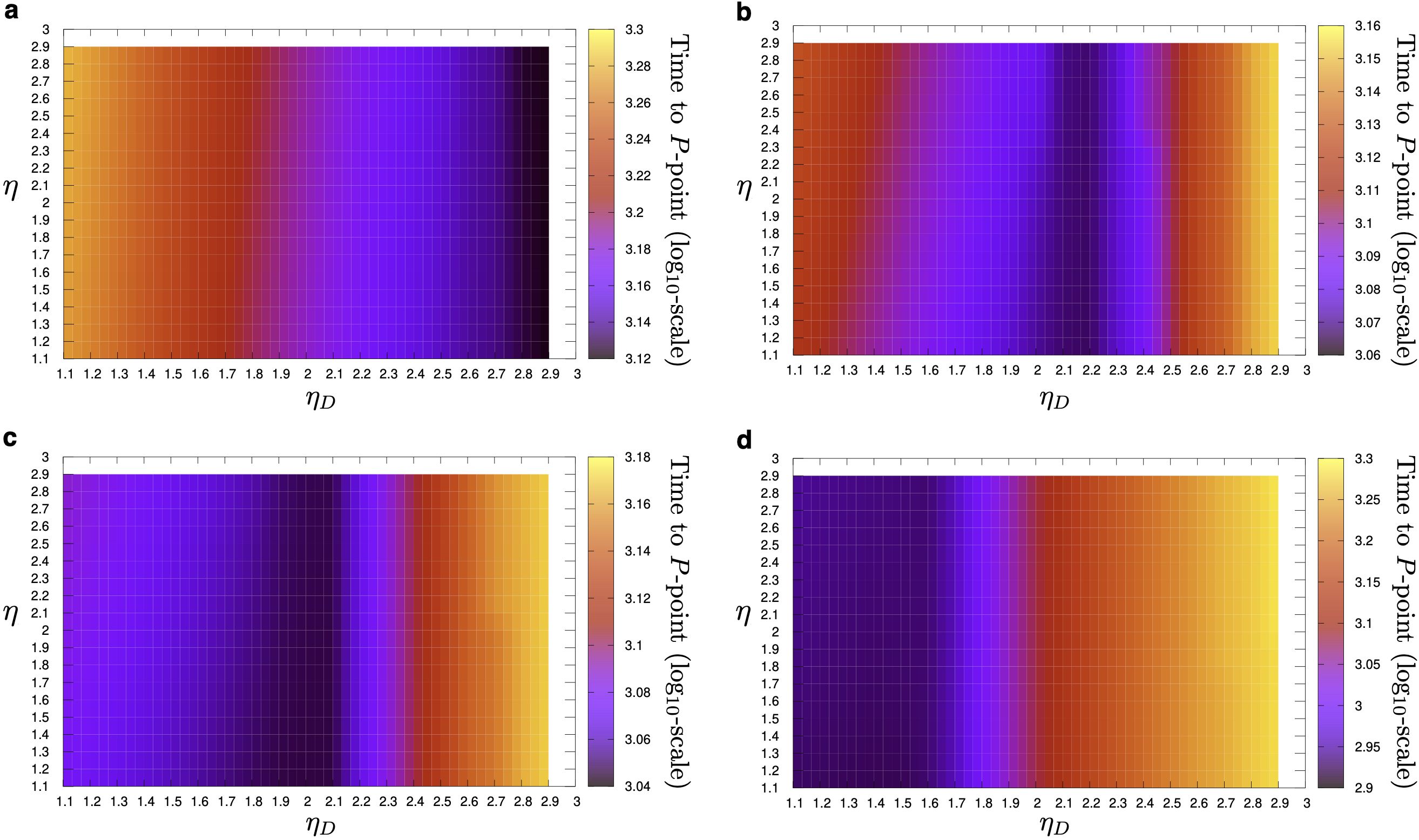}
\caption{Transient times (in $\log_{10}$-scale) to reach a distance $10^{-10}$ from the (attracting and unique) $P$-point by the orbit with initial conditions $(V(0),S(0),D(0))=(0.1, 0.05, 0)$ computed in the $(\eta_D,\eta)$ plane. (a) Times obtained from replication rates from Section~\ref{values:HBV-HDV}, given by $\alpha=1$, $\beta=1.84$ and $\gamma=2.40$. Other values for $\gamma$ are shown in panels (b) $3.5$; (c) $4$, and (d) $10$. In all of the panels we have used $\mu=0.47$ and $\eps=3 \cdot 10^{-2}$.}
\label{fig:times:HBV}
\end{figure}

\subsubsection{Time to equilibria: case $\gamma>\beta$.} \label{values:HBV-HDV}

As discussed in the Introduction, the most common situation is that DIs  take advantage of their shorter genome to replicate faster than their parental HV. Furthermore, their genomes are shorter than the ones characteristic of the large linear satellites (median length $1.1 \pm 0.4$ Kb ($\pm$ IQR)), but not necessarily so for the small linear (median length $0.4 \pm 0.1$ Kb) nor the virusoid (median length $0.3 \pm 0.1$ Kb) satellites \cite{Badar2021} accompanying the HV. Mathematically, this translates into the fact that $\gamma \gtrsim \beta$. This implies that full coexistence i.e., $Q$-points, would rarely exist, and equilibrium scenarios with no satRNA would be the most common outcome. That is, DIs are efficient enough to outcompete the longer satRNA from this steady state solution. Despite this equilibrium situation, satRNAs could persist in a transitory way in the system for very long times, as mentioned above. To illustrate these transients, we have chosen the clinically relevant case of hepatitis B virus (HBV), its defective D-RNAs, and the satellite hepatitis delta virus (HDV). In order to numerically simulate this case, it is reasonable to assume that their replication rates are proportionally inverse to their genome's length. If they are
\[
\textrm{HBV}: [3017,3248], \qquad \textrm{HDV}: [1679,1682], \qquad
\textrm{D-RNA}: \sim 1290
\]
nucleotides long, we can obtain the following estimates for the three replication rates:
\[
\alpha=1, \qquad \beta \simeq \frac{3100}{1680} \simeq 1.84, \qquad
\gamma \simeq \frac{3100}{1290} \simeq 2.4.
\]
Figure~\ref{fig:times:HBV} displays transient times for this particular case where $\gamma > \beta$. Specifically, panel (a) shows the result for the replication rates listed above for the HBV-HDV case while in the other panels $\gamma$ has been further increased. In all cases the orbit starts with initial conditions $(V(0),S(0),D(0))=(0.1, 0.05, 0)$ and the distance from the (attracting) $P$-point at which numerical integration stops is $10^{-10}$. Notice that, although the coordinates of the $P$-point do not depend on $\eta$, the vector field does. This produces slightly variations in these transient times. Moreover, the results show that $\eta_D$ has a stronger impact on these transient times towards the outcompetition of the satellite in comparison to $\eta$. Observe also that these transients  vary in a nontrivial way as $\gamma$ increases. For instance, for the values of $\gamma$ and $\beta$ chosen as representatives of the HBV-HDV virus system, the longest times are found at low interference values of DI RNAs. When DIs replication is further increased, this region with longer transients moves to larger $\eta_D$ values. This result is somehow counter-intuitive since faster DIs' replication should involve faster satellites extinction, but this outcome is probably counteracted by the larger interference of DIs on the HV. Regardless of these results, we must notice that the difference in the length of the transients between the yellow-orange and the black-blue scales is not very large (as compared to those of Fig.~\ref{fig:timetoQpoints}). Although the replication cycle of HBV and HDV is rather more complex and our model may be very limited due to its level of abstraction, it indicates that the faster replicator (in this case the DIs) will typically outperform the other sub-viral element, so the dynamics follow Gause's competitive exclusion principle.

\begin{comment}
#define ALPHA  1.             // PMV: 4326.       // HBV  [3017,3248] --> 3100.
#define BETA   3100./1680.    //      4326./824.  //      [1679,1682] --> 3100./1680.
#define GAMMA  10.            //      4326./400.  //      1290        --> 3100./1290.
#define EPSILON 3.e-2

#define ETA_MIN 1.1
#define ETA_MAX 3.
#define ETAD_MIN 1.1
#define ETAD_MAX 3.

#define ACCURACY_FIXED_POINT 1.e-13

// for the plot (window)
#define MU  0.47
#define GRID_STEP 1.e-1

// initial conditions for the orbit selected
#define V_INITIAL  0.1   // 0.01
#define S_INITIAL  0.05   // 0.05
#define D_INITIAL  0.

#define DISTANCE_TO_EQPOINT 1.e-10
\end{comment}

\section{Conclusions}
DI RNAs are an unavoidable consequence of the error-prone replication of RNA viruses and retroviruses.  The impact of these defectors in the population dynamics of their parental virus have been deeply studied~\cite{Szathmary1993,BanghamKirkwood1990,Kirkwood1994,Sardanyes2010,Zwart2013,Chao2017}. However, viral infections, specially in plants, are more complex and contain additional genetic elements that are unrelated with the virus: the satellite RNAs (satRNAs) and the satellite viruses. Both DIs and satellites share a common feature, they need the wild type virus for their own replication since they lack essential genes such as those coding for the viral polymerases or for the coat proteins. Hence, they need to co-infect with the wild type virus (helper virus, HV) to complete their replication/infection cycle. The presence of these extranumerary elements have been shown to deeply affect the virulence of infection \cite{Simon2004,Gnanasekaran2019,Taylor2020}, in some cases exacerbating symptoms while in other resulting in their attenuation. Therefore, the interaction between satellites and their HV ranges from commensalism to parasitism.  Here, we present a simple, yet dynamically rich, model of infection with a HV, a generic satellite and the DI RNAs generated from the HV.  All three RNA species compete for limited host resources, thus we implicitly assume the satRNA acts as a hyperparasite.

Analytical and numerical explorations of the model show three possible stable states: (i) full extinction; (ii) outcompetition of the satRNA by the duo HV-DI RNAs; and (iii) coexistence of the three replicators. A rather small region of bistability involving coexistence of states (ii) and (iii) has been found, having the fixed point responsible for scenario (ii) a very small basin of attraction. 
We have analytically found the condition under which the three replicators can go to extinction, showing that there is a critical rate of DIs production that only depends on the balance between the degradation and replication of the HV through the equation $\mu_c = 1 - \eps / \alpha$. This means that when the rate of production of DIs overcomes the critical condition $\mu_c$, a full virus clearance occurs through a transcritical bifurcation. Note that this critical value does not depend at all on the satRNAs parameters. We have also identified that the majority of transitions between scenarios (i) to (iii) are given by transcritical bifurcations, except for the tiny bistability region, where saddle-node bifurcations are found. Most remarkably from an applied perspective, we found conditions in which the WT virus takes advantage of the unavoidable production of DIs to outcompete the satellite.  Indeed, the strength of this outcompetition effect becomes stronger as the difference in lenght between the DI RNAs and the satRNAs increases: large linear satellites, and the special case of the HDV virusoid, outcompetition takes place rapidly ($\beta < \gamma$), while for small linear satellites, it might take longer ($\beta \approx \gamma$).

The model here presented is a minimum one that lacks of mechanistic details and only focuses in replication interactions. An obvious extension of the study of these interesting multi-species viral models would consist in including proteins in the picture: HV will encode for replication and encapsidation machinery whereas DI RNAs would simply kidnap these proteins for their own replication and encapsidation. In such mechanistic model, satellite viruses and satRNAs would not be collapsed into a single category, as we have done here, but represented by two different molecular species, one encoding for some protein (satellite virus) and other do not encoding for any factor (satRNAs).  Another possible extension of our model would consist in imposing a second layer of complexity involving eco-evolutionary dynamics, \emph{e.g.} a multi-strain SIR-like model.

Hyperparasitism could potentially play a key role in biological control of viral infections \cite{Sandhu2021} by reducing the deleterious impact of the WT virus in its host, also hampering its transmission.  Indeed, this principle is the ground for the recent development of antiviral therapies based in the generation of engineered artificial DI RNAs that strongly interfere with the target virus \cite{Notton2014,Tanner2016}.  These novel approaches, however, need of additional careful theoretical considerations, as recent eco-evolutionary models have shown that introducing a hyperparasite into the original host-parasite system results in a shift of the evolutionarily optimal virulence of the pathogen toward higher values \cite{Sandhu2021}.  Our results suggest that, in the case of highly mutable RNA viruses, the constitutive production of DI RNAs may contribute to avoid the establishment of a hyperparasite competing for helper virus resources.

\section*{Acknowledgements}
JTL has been funded by Agencia Estatal de Investigaci\'on (AEI) grants PGC2018-098676-B-100 and PID2021-122954NB-I00 funded by MCIN/AEI/10.13039/501100011033/ and "ERDF A way of making Europe", and by "Ayudas para la Recualificaci\'on del Sistema Universitario Espa\~nol 2021-2023". JS has been funded by a Ram\'on y Cajal Fellowship (RYC-2017-22243) funded by MCIN/AEI/10.13039/501100011033 "FSE invests in your future". JS and TA also acknowledge grant RTI2018-098322-B-I00 funded by MCIN/AEI/10.13039/501100011033/ and "ERDF A way of making Europe", and Generalitat de Catalunya CERCA Program for institutional support. SFE has been supported by AEI-FEDER grant PID2019-103998GB-I00 and Generalitat Valenciana grant PROMETEO/2019/012. We also acknowledge support from Mar\'ia de Maeztu Program for Units of Excellence in R\&D grant CEX2020-001084-M (JTL, JS, TA). JTL also thanks the hospitality of Laboratorio Subterraneo de Canfranc and I$^2$SysBio as hosting institutions of this grant. We want to thank Lluís Alsedà, Julia Hillung, Axel Masó, Juan C. Muñoz, María J. Olmo, Josep Quer and members of his group, Pau Reig, Francisco Rodríguez Frías, David Romero, and Marco Vignuzzi for interesting and stimulating discussions about defective viral genomes.

\end{document}